\documentclass[pdflatex,sn-mathphys-num]{sn-jnl}
\usepackage{amsmath,amssymb,mathtools,bm,microtype}
\usepackage{xcolor}
\newcommand{\rev}[1]{#1}

\theoremstyle{thmstyleone}
\newtheorem{theorem}{Theorem}[section]
\newtheorem{lemma}[theorem]{Lemma}

\newtheorem{corollary}[theorem]{Corollary}
\theoremstyle{thmstyletwo}

\numberwithin{equation}{section}
\begin{document}
\title[Capacitary Hardy inequalities]{Capacitary Hardy inequalities, Dirichlet boundary conditions, and the limit as $p\downarrow1$}
\author*[1]{\fnm{Takanobu} \sur{Hara}}\email{takanobu.hara.math@gmail.com}
\affil*[1]{\orgdiv{Mathematical Institute}, \orgname{Tohoku University}, \orgaddress{\street{6-3, Aramaki Aza-Aoba, Aoba-ku}, \city{Sendai}, \state{Miyagi}, \postcode{980-8578}, \country{Japan}}}
\date{November 12, 2026}
\abstract{
We use nonlinear Dirichlet resolvents to prove Hardy inequalities with capacity-based and resolvent-based weights for $1<p<\infty$ on arbitrary open sets. On bounded Sobolev extension domains, integrability against either weight characterizes Sobolev functions satisfying homogeneous Dirichlet conditions on arbitrary closed subsets of the boundary. On bounded Lipschitz domains, the resolvent-based weight $\bm{W}_{D,p}$ converges locally in $L^1$ to the negative of Barozzi's interior variational mean curvature as $p\downarrow1$.
}
\maketitle

\section{Introduction}\label{sec:intro}
\begingroup\color{black}

The validity of Hardy inequalities in higher-dimensional Euclidean spaces is closely related to potential theory.
Classical proofs of the distance Hardy inequality use strong barriers \cite{Ancona,HaraBarrier}, 
self-improvement of uniform fatness \cite{Lewis,Mikkonen,KM}, or boundary coverings \cite{Wannebo}. 
Uniform fatness characterizes the pointwise Hardy inequality \cite[Theorem~1.1]{KLT}.
For an overview, see \cite{KLV}.

Chen, Xiao, Yang, Yuan, and Zhang introduced a probabilistic approach based on the absorption time of killed Brownian motion.
They used this approach to solve Maz'ya's open problem \cite[Section~3.8, Problem~8]{Mazya} by proving the capacitary-distance Hardy inequality on arbitrary open sets $D\subset\mathbb R^n$ for $n\ge3$ \cite[Theorem~1.3]{CXYYZ}:
\begin{equation}\label{eq:cxintro}
 \int_D\frac{|u(x)|^2}{d_\alpha(x)^2}\,dx
 \le\frac{C_n}{\alpha^2}\int_D|\nabla u|^2\,dx,
 \qquad u\in C_c^\infty(D),\quad 0<\alpha\le1.
\end{equation}
Here $d_\alpha(x)$ is the infimum of radii for which the Newtonian capacity of $\overline{B(x,r)\setminus D}$ is at least $\alpha$ times that of $B(x,r)$. 
They also proved that the dependence on $\alpha^{-2}$ is sharp. Their variable-time semigroup estimate gives a Hardy inequality weighted by the reciprocal of an absorption-time quantile for killed Brownian motion. They control that quantile by showing that a capacity ratio of at least $\alpha$ forces absorption with probability at least $c_n\alpha$ by time $C_n r^2$ \cite[Theorems~1.1 and~1.2]{CXYYZ}.

Inspired by this result, the squared relative capacity ratio was integrated over scales in \cite[Proposition~8.4 and Appendix~A]{Hara} to obtain a Hardy weight and establish
\[
 H_F^1(D)=\left\{u\in H^1(D):\int_D |u(x)|^2 \rev{\bm{V}}_{F,2}(x)\,dx<\infty\right\}
\]
for bounded $H^1$-extension domains in dimension $n\ge2$ and closed $F\subset\partial D$, with the weight expressed in terms of the equivalent capacity density in \eqref{eq:theta}. Here $H_F^1(D)$ is the closure in $H^1(D)$ of the restrictions to $D$ of functions in $C_c^\infty(\mathbb R^n\setminus F)$, as in \eqref{eq:partialspace} with $p=2$. For comparison, the criterion of Kinnunen and Martio \cite[Theorem~3.13]{KM} gives $u\in W^{1,p}_0(D)$ for $u\in W^{1,p}(D)$ if $u/\operatorname{dist}(\cdot,\partial D)$ belongs to weak $L^p(D)$.

For $n\ge3$,
\[
 \int_0^1\frac{2\alpha}{d_\alpha(x)^2}\,d\alpha\asymp_n \rev{\bm{V}}_{\mathbb R^n\setminus D,2}(x).
\]
By the layer-cake formula, this equals $2\int_0^\infty(\sup_{0<s<r}q_x(s))^2\,dr/r^3$, where $q_x(r)$ is the Newtonian capacity ratio defining $d_\alpha(x)$. Monotonicity of the numerator bounds this squared supremum between $q_x(r)^2$ and $C_n\int_0^{2r}q_x(t)^2\,dt/t$ for almost every $r$. Tonelli's theorem and condenser comparison give the equivalence. Restricting the threshold integral to $(0,\alpha)$ recovers the upper bound in \eqref{eq:cxintro}. Integrating that bound against $2\alpha\,d\alpha$ gives a logarithmic divergence.

We extend the capacity-integrated Hardy inequality and partial zero-trace characterization to $1<p<\infty$. The absorption-time argument bounds the measure defined by the reciprocal squared absorption-time scale in terms of relative capacity \cite[(4.5)]{CXYYZ}. For $p=2$, the corresponding estimate on bounded open sets was obtained by integrating a linear resolvent identity over scales \cite[Lemma~A.2]{Hara}. Its nonlinear extension gives a capacitary estimate for $\rev{\bm{W}}_{D,p}$, and a local estimate transfers the Hardy bound to $\rev{\bm{V}}_{\mathbb R^n\setminus D,p}$. The weights are comparable for every $1<p<\infty$ by Theorem~\ref{thm:comparison}.

Beyond establishing two Hardy inequalities, Theorem~\ref{thm:comparison} provides a pointwise two-sided comparison between the two weights on arbitrary open sets, without any boundary regularity assumptions and including the critical case $n=p=2$. In the linear case $p=2$, let $\tau_D^{(1)}$ and $\tau_D^{(2)}$ be independent exit times of Brownian motions started at $x$, each with generator $\Delta$. Since $z_{r,2}^D(x)=\mathbb E_x[e^{-\tau_D/r^2}]$, the substitution $t=r^{-2}$ yields
\begin{equation}\label{eq:introexit}
 \bm{W}_{D,2}(x)=\mathbb E_x\!\left[\frac{1}{\tau_D^{(1)}+\tau_D^{(2)}}\right]\asymp_n\bm{V}_{\mathbb R^n\setminus D,2}(x),
\end{equation}
where $1/\infty=0$. Thus the expected reciprocal of the sum of two independent exit times is comparable to the squared capacity density of the complement integrated over all scales. For general $p$, the main difficulty is the reverse estimate $\bm{W}_{D,p}\le C_{n,p}\bm{V}_{\mathbb R^n\setminus D,p}$, since the linear annular decomposition of equilibrium potentials is unavailable. Section~\ref{sec:comparison} establishes this estimate through an exponentially weighted bound for the nonlinear equilibrium measure and convolution estimates on logarithmic scales.

On a bounded $W^{1,p}$-extension domain, Theorem~\ref{thm:trace} gives, for closed $F\subset\partial D$,
\begin{equation}\label{eq:introtrace}
 \begin{aligned}
 W_F^{1,p}(D)
 &=\{u\in W^{1,p}(D): |u|^p\rev{\bm{V}}_{F,p}\in L^1(D)\}\\
 &=\{u\in W^{1,p}(D): |u|^p\rev{\bm{W}}_{\mathbb R^n\setminus F,p}\in L^1(D)\}.
 \end{aligned}
\end{equation}
The constants in these weighted estimates are independent of $F$.

Zero-trace characterizations by weighted integrability go back to Kadlec and
Kufner \cite{KadlecKufnerI,KadlecKufnerII}. Later results include weak-$L^p$
and Lorentz-space criteria \cite{KM,NekvindaTurcinova}, as well as
higher-order analogues \cite{EdmundsNekvinda}. For conditions imposed on a
closed part of the boundary, previous work includes distance Hardy inequalities
under geometric assumptions \cite{EgertHallerDintelmannRehberg} and zero-trace
characterizations involving distance weights under additional regularity
assumptions \cite{EgertTolksdorf}. In comparison,
\eqref{eq:introtrace} uses capacity-based weights and permits any closed
$F\subset\partial D$ on a bounded Sobolev extension domain.

On bounded Lipschitz domains, $\rev{\bm{W}}_{D,p}$ converges locally in $L^1$ as $p\downarrow1$ to $-\rev{\mathsf{H}}_D$, where $\rev{\mathsf{H}}_D$ is Barozzi's interior variational mean curvature \cite{Barozzi}. Barozzi's perimeter bound for $-\rev{\mathsf{H}}_D$ yields a sharp BV Hardy inequality.
\endgroup

\section{Nonlinear resolvents and Hardy inequalities}\label{sec:finite}
Let $D\subsetneq\mathbb R^n$ be nonempty and open, where $n\ge2$. Set $\delta_D(x)=\operatorname{dist}(x,\mathbb R^n\setminus D)$. Throughout this section, $1<p<\infty$. The space $W^{1,p}_0(D)$ is the closure of $C_c^\infty(D)$ in $W^{1,p}(D)$. \begingroup\color{black}
For $r>0$ and compact $K\subset\mathbb R^n$, define the whole-space Sobolev capacity at scale $r$ by
\[
 \operatorname{Cap}_{p,r}(K)=
 \inf_{\substack{\varphi\in C_c^\infty(\mathbb R^n)\\
 \varphi\ge1\text{ near }K}}
 \int_{\mathbb R^n}\bigl(|\nabla\varphi|^p+r^{-p}|\varphi|^p\bigr)\,dy.
\]
Extend it by outer regularization. A property holds $p$-quasieverywhere if it holds outside a set of $\operatorname{Cap}_{p,1}$-capacity zero. Scaling and translation invariance give
\begin{equation}\label{eq:ballcap}
 \operatorname{Cap}_{p,r}(B(x,r))=r^{n-p}\operatorname{Cap}_{p,1}(B(0,1)).
\end{equation}
For a closed set $F\subset\mathbb R^n$, define the Sobolev capacity density and its Hardy weight by
\begin{equation}\label{eq:theta}
 \begin{aligned}
 \Theta_{F,p}(x,r)&=\frac{\operatorname{Cap}_{p,r}(F\cap B(x,r))}
 {\operatorname{Cap}_{p,r}(B(x,r))},\\
 \rev{\bm{V}}_{F,p}(x)&=p\int_0^\infty\Theta_{F,p}(x,r)^{p/(p-1)}\frac{dr}{r^{p+1}}.
 \end{aligned}
\end{equation}
If $K\subset B(x,r)$, then zero extension, the Poincar\'e inequality, and a cutoff yield
\begin{equation}\label{eq:relative-capacity}
 \operatorname{Cap}_{p,r}(K)\asymp_{n,p}
 \operatorname{cap}_p(K,B(x,2r)),
\end{equation}
where $\operatorname{cap}_p(K,U)$ is the infimum of $\int_U|\nabla\varphi|^p\,dy$ over $\varphi\in C_c^\infty(U)$ with $\varphi\ge1$ near $K$.
\endgroup
Since $\Theta_{\mathbb R^n\setminus D,p}(x,r)=0$ for $r<\delta_D(x)$ and $0\le\Theta_{\mathbb R^n\setminus D,p}\le1$, we have
\begin{equation}\label{eq:capupper}
 0\le \rev{\bm{V}}_{\mathbb R^n\setminus D,p}(x)\le\delta_D(x)^{-p}.
\end{equation}

First suppose that $D$ is bounded. For $r>0$ and $f\in L^p(D)$, let $J_{r,p}f$ be the unique minimizer of
\begin{equation}\label{eq:resolvent}
 \int_D|\nabla v|^p\,dx+r^{-p}\int_D|f-v|^p\,dx,
 \qquad v\in W^{1,p}_0(D).
\end{equation}
\rev{Set $z_{r,p}^D=1-J_{r,p}1$.} Truncation and the comparison principle give $0\le z_{r,p}^D\le1$. Its continuous interior representative satisfies
\begin{equation}\label{eq:absorption}
 -\Delta_p z_{r,p}^D+r^{-p}(z_{r,p}^D)^{p-1}=0\quad\text{in }D,
 \qquad 1-z_{r,p}^D\in W^{1,p}_0(D).
\end{equation}
For an unbounded $D$, take $D_j=D\cap B(0,j)$ and define $z_{r,p}^D$ as the decreasing limit of $z_{r,p}^{D_j}$ at interior points. Comparison shows that this definition is independent of the bounded exhaustion. Write
\begin{equation}\label{eq:massive}
 \rev{\bm{W}}_{D,p}(x)=p\int_0^\infty z_{r,p}^D(x)^p\frac{dr}{r^{p+1}}.
\end{equation}
If $D_1\subset D_2$, then $z_{r,p}^{D_2}\le z_{r,p}^{D_1}$ on $D_1$.

\begin{theorem}\label{thm:hardy}
Let $D\subsetneq\mathbb R^n$ be open and let $1<p<\infty$. For every $u\in W^{1,p}_0(D)$,
\begin{align}
 \int_D \rev{\bm{W}}_{D,p}|u|^p\,dx&\le 2^p(2^p-1)\int_D|\nabla u|^p\,dx,\label{eq:massivehardy}\\
 \int_D \rev{\bm{V}}_{\mathbb R^n\setminus D,p}|u|^p\,dx&\le C_{n,p}\int_D|\nabla u|^p\,dx.\label{eq:caphardy}
\end{align}
Moreover, $\rev{\bm{V}}_{\mathbb R^n\setminus D,p}(x)\le C_{n,p}\rev{\bm{W}}_{D,p}(x)$ for every $x\in D$.
\end{theorem}

\begin{lemma}\label{lem:identity}
Suppose that $D$ is bounded and $f\in W^{1,p}_0(D)$. Then
\begin{equation}\label{eq:identity}
 p\int_0^\infty\|f-J_{r,p}f\|_{L^p(D)}^p\frac{dr}{r^{p+1}}
 =\int_D|\nabla f|^p\,dx.
\end{equation}
\end{lemma}
\begin{proof}
\begingroup\color{black}
For $\lambda>0$, set
\[
 \Phi(\lambda)=\inf_{v\in W^{1,p}_0(D)}
 \left\{\int_D|\nabla v|^p\,dx+\lambda\int_D|f-v|^p\,dx\right\}.
\]
The function $\Phi$ is increasing and concave. Comparing minimizers at two parameters gives
$\Phi'(\lambda)=\|f-J_{\lambda^{-1/p},p}f\|_{L^p(D)}^p$ for almost every $\lambda>0$.
The competitor $v=0$ gives $\Phi(0+)=0$. The competitor $v=f$ bounds $\Phi(\lambda)$ by the energy of $f$. As $\lambda\to\infty$, the minimizers converge to $f$ in $L^p(D)$ and have bounded energies. Weak compactness and lower semicontinuity therefore give
$\Phi(\infty)=\int_D|\nabla f|^p\,dx$.
Integrating $\Phi'$ and substituting $\lambda=r^{-p}$ proves \eqref{eq:identity}.
\endgroup
\end{proof}

\begin{proof}[Proof of \eqref{eq:massivehardy}]
\begingroup\color{black}
First let $D$ be bounded. For compact $K\subset D$, choose $f\in W^{1,p}_0(D)$ with $0\le f\le1$ and $f=1$ in a neighborhood of $K$. Order preservation gives $J_{r,p}f\le J_{r,p}1$. Hence $f-J_{r,p}f\ge z_{r,p}^D$ on $K$. Integrating and applying Lemma~\ref{lem:identity} gives
\begin{equation}\label{eq:weighted-capacity-bound}
 \int_K \rev{\bm{W}}_{D,p}\,dx\le\int_D|\nabla f|^p\,dx.
\end{equation}
For unbounded $D$, apply \eqref{eq:weighted-capacity-bound} on $D\cap B(0,j)$ to each compactly supported competitor and pass to the decreasing limit of the responses.

For $u\in C_c^\infty(D)$ and $j\in\mathbb Z$, put
\[
 L_j=\{|u|\ge2^j\},\qquad
 h_j=\min\{1,2^{1-j}(|u|-2^{j-1})_+\}.
\]
Apply \eqref{eq:weighted-capacity-bound} with $K=L_j$ and $f=h_j$. To justify this application, first truncate $(1+\varepsilon)h_j$ at $1$, approximate the result by smooth admissible functions, and then let $\varepsilon\downarrow0$. Hence $\int_{L_j}\rev{\bm{W}}_{D,p}\,dx\le\int_D|\nabla h_j|^p\,dx$.
\endgroup
The dyadic sums satisfy
\[
 \sum_{j\in\mathbb Z}2^{jp}\rev{\bm{1}}_{L_j}\ge\frac{|u|^p}{2^p-1},
 \qquad
 \sum_{j\in\mathbb Z}2^{jp}|\nabla h_j|^p\le2^p|\nabla u|^p
\]
almost everywhere. Multiplication by $2^{jp}$, summation, and Tonelli's theorem prove \eqref{eq:massivehardy} for smooth $u$.
Approximation in $W^{1,p}_0(D)$ and Fatou's lemma complete the proof.
\end{proof}

\begin{lemma}[Capacitary barrier]\label{lem:barrier}
For every $x\in D$ and $r>0$,
\begin{equation}\label{eq:barrier}
 z_{r,p}^D(x)\ge c_{n,p}\Theta_{\mathbb R^n\setminus D,p}(x,r)^{1/(p-1)}.
\end{equation}
\end{lemma}
\begin{proof}
\begingroup\color{black}
Fix a compact set $K\subset (\mathbb R^n\setminus D)\cap B(x,r)$. Let $w\in W^{1,p}_0(B(x,8r))$ minimize
\[
 \int_{B(x,8r)}\bigl(|\nabla w|^p+r^{-p}|w|^p\bigr)\,dy
\]
subject to $w\ge1$ $p$-quasieverywhere on $K$. Truncation gives $0\le w\le1$. The function $w$ is a supersolution of the absorption equation in the ball and a solution outside $K$.

After scaling to $r=1$, the weak Harnack inequality \cite[Theorem~1.2]{Trudinger} with exponent $p-1$ gives \eqref{eq:localmass}. The Harnack inequality on balls away from $K$ \cite[Theorem~1.1]{Trudinger} and the Caccioppoli inequality then give \eqref{eq:localgradient}. In the original scale,
\begin{align}
 \int_{B(x,4r)}w^{p-1}\,dy&\le Cr^nw(x)^{p-1},\label{eq:localmass}\\
 \int_{B(x,3r)\setminus B(x,2r)}|\nabla w|^p\,dy&\le Cr^{n-p}w(x)^p.\label{eq:localgradient}
\end{align}
Here $C$ depends only on $n$ and $p$.

Choose $\eta\in C_c^\infty(B(x,3r))$ with $0\le\eta\le1$, $\eta=1$ on $B(x,2r)$, and $|\nabla\eta|\le C/r$. The admissible variation $\eta^p(1-w)$ vanishes on $K$. The variational equation, H\"older's inequality, and the estimates \eqref{eq:localmass} and \eqref{eq:localgradient} give
\[
 \begin{aligned}
 \int_{B(x,3r)}\eta^p|\nabla w|^p\,dy
 &\le p\int_{B(x,3r)}\eta^{p-1}|\nabla w|^{p-1}|\nabla\eta|\,dy
       +r^{-p}\int_{B(x,3r)}w^{p-1}\,dy\\
 &\le Cr^{n-p}w(x)^{p-1}.
 \end{aligned}
\]
The zero extension of $\eta w$ is admissible for $\operatorname{Cap}_{p,r}(K)$. The estimates \eqref{eq:localmass} and \eqref{eq:localgradient}, together with $w^p\le w^{p-1}$ and $w(x)\le1$, give $\operatorname{Cap}_{p,r}(K)\le Cr^{n-p}w(x)^{p-1}$.
Comparison on $D\cap B(x,8r)$ gives $w\le z_{r,p}^D$. Taking compact approximations to $(\mathbb R^n\setminus D)\cap B(x,r)$ and using \eqref{eq:ballcap} proves \eqref{eq:barrier}. For an unbounded $D$, apply comparison in the bounded exhaustion and pass to its decreasing limit.
\endgroup
\end{proof}

\begin{proof}[Proof of \eqref{eq:caphardy}]
\rev{Integrating \eqref{eq:barrier} to the power $p$ gives $\rev{\bm{V}}_{\mathbb R^n\setminus D,p}\le C_{n,p}\rev{\bm{W}}_{D,p}$. Together with \eqref{eq:massivehardy}, this proves \eqref{eq:caphardy}.}
\end{proof}

\begin{corollary}[Distance Hardy inequality]\label{cor:distance-hardy}
Let $D\subsetneq\mathbb R^n$ be nonempty and open, and let $1<p<\infty$.
Suppose either that $p>n$, or that there is $\eta>0$ such that
\[
 \operatorname{cap}_p((\mathbb R^n\setminus D)\cap B(y,s),B(y,2s))
 \ge\eta\operatorname{cap}_p(B(y,s),B(y,2s))
\]
for every $y\in\partial D$ and $s>0$. Then
\begin{equation}\label{eq:distancehardy}
 \int_D\frac{|u(x)|^p}{\delta_D(x)^p}\,dx
 \le C\int_D|\nabla u|^p\,dx,
 \qquad u\in W_0^{1,p}(D),
\end{equation}
where $C$ depends only on $n,p$ if $p>n$, and on $n,p,\eta$ under the capacity density assumption.
\end{corollary}
\begin{proof}
If $p>n$, the positivity of the $p$-capacity of a singleton, scaling, and \eqref{eq:ballcap} give
$\Theta_{\mathbb R^n\setminus D,p}(x,r)\ge c_{n,p}$ for $r>\delta_D(x)$.
Under the capacity density assumption, choose $y\in\partial D$ with $|x-y|=\delta_D(x)$.
For $r\ge2\delta_D(x)$, we have $B(y,r/2)\subset B(x,r)$.
The capacity density condition, \eqref{eq:relative-capacity}, and \eqref{eq:ballcap} give
$\Theta_{\mathbb R^n\setminus D,p}(x,r)\ge c_{n,p}\eta$.
In either case, integration in \eqref{eq:theta} gives
$\bm{V}_{\mathbb R^n\setminus D,p}(x)\ge c\delta_D(x)^{-p}$.
The assertion follows from \eqref{eq:caphardy}; the case $p>n$ is classical \cite[Section~3]{KM}.
\end{proof}

\begingroup\color{black}
\section{Capacitary characterization of Dirichlet boundary conditions}\label{sec:trace}
Let $D\subset\mathbb R^n$ be a bounded $W^{1,p}$-extension domain, where $1<p<\infty$. Fix a bounded linear extension operator
$\mathcal E_D:W^{1,p}(D)\to W^{1,p}(\mathbb R^n)$.
For a closed set $F\subset\partial D$, define
\begin{equation}\label{eq:partialspace}
 W_F^{1,p}(D)=\overline{\{\varphi|_D:\varphi\in C_c^\infty(\mathbb R^n\setminus F)\}}^{\,W^{1,p}(D)}.
\end{equation}
We use the weight $\rev{\bm{V}}_{F,p}$ defined in \eqref{eq:theta}. The resolvent-based weight is $\rev{\bm{W}}_{\mathbb R^n\setminus F,p}$. For $F=\varnothing$, set $\rev{\bm{W}}_{\mathbb R^n,p}=0$.

\begin{theorem}[Characterization by Hardy weights]\label{thm:trace}
Under these assumptions, the following conditions are equivalent for $u\in W^{1,p}(D)$.
\begin{enumerate}
\item $u\in W_F^{1,p}(D)$.
\item The quasicontinuous representative of $\mathcal E_Du$ vanishes $p$-quasieverywhere on $F$.
\item $\int_D|u|^p\rev{\bm{V}}_{F,p}\,dx<\infty$.
\item $\int_D|u|^p\rev{\bm{W}}_{\mathbb R^n\setminus F,p}\,dx<\infty$.
\end{enumerate}
For $u\in W_F^{1,p}(D)$,
\begin{equation}\label{eq:partialbound}
 \int_D|u|^p(\rev{\bm{V}}_{F,p}+\rev{\bm{W}}_{\mathbb R^n\setminus F,p})\,dx
 \le C\|u\|_{W^{1,p}(D)}^p,
\end{equation}
where $C$ depends on $n$, $p$, $D$, and the norm of $\mathcal E_D$, and is independent of $F$.
\end{theorem}

We first record the boundary-value properties needed in the proof. The extension property implies the measure density estimate \cite[Theorem~2]{HKT}
\begin{equation}\label{eq:densitytrace}
 |D\cap B(y,r)|\ge br^n\qquad(y\in\partial D,\ 0<r<r_D)
\end{equation}
for some $b>0$ and $r_D>0$. For $U\in W^{1,p}(\mathbb R^n)$, its quasicontinuous representative $\widetilde U$ satisfies the capacitary Lebesgue-point property \cite[Theorem~4.5 and the preceding paragraph]{KL}
\begin{equation}\label{eq:lebesguetrace}
 \lim_{r\downarrow0}\frac1{|B(y,r)|}\int_{B(y,r)}|U-\widetilde U(y)|\,dx=0
 \quad\text{for }p\text{-quasievery }y.
\end{equation}
If two extensions agree on $D$, \eqref{eq:densitytrace} and \eqref{eq:lebesguetrace} imply that their representatives agree $p$-quasieverywhere on $\partial D$. Thus condition (2) is independent of the chosen extension.

\begin{proof}[Proof of (1)$\Leftrightarrow$(2) and \eqref{eq:partialbound}]
The case $F=\varnothing$ follows from the extension property and smooth approximation, so assume $F\ne\varnothing$. Suppose (1) holds. Choose approximating restrictions $\varphi_j|_D$ in \eqref{eq:partialspace}. The extensions $\mathcal E_D(\varphi_j|_D)$ converge to $\mathcal E_Du$ in $W^{1,p}(\mathbb R^n)$. Their representatives vanish on $F$ outside sets of zero capacity by boundary-value uniqueness. A subsequence of these extensions converges to $\mathcal E_Du$ $p$-quasieverywhere, proving (2).

Conversely, suppose (2) holds. Choose $\chi\in C_c^\infty(\mathbb R^n)$ equal to $1$ near $\overline D$. Put $g=\chi\mathcal E_Du$. The quasieverywhere characterization of $W^{1,p}_0$ \cite[Theorem~4.5]{HKM} gives $g\in W^{1,p}_0(\mathbb R^n\setminus F)$. Approximation in this space by $C_c^\infty(\mathbb R^n\setminus F)$ and restriction to $D$ prove (1). Theorem~\ref{thm:hardy}, applied on $\mathbb R^n\setminus F$, also gives
\[
 \int_D|u|^p(\rev{\bm{V}}_{F,p}+\rev{\bm{W}}_{\mathbb R^n\setminus F,p})\,dx
 \le C_{n,p}\int_{\mathbb R^n}|\nabla g|^p\,dx
 \le C\|u\|_{W^{1,p}(D)}^p.
\]
This proves \eqref{eq:partialbound} and the implications (1)$\Rightarrow$(3) and (1)$\Rightarrow$(4). The comparison $\rev{\bm{V}}_{F,p}\le C_{n,p}\rev{\bm{W}}_{\mathbb R^n\setminus F,p}$ proves (4)$\Rightarrow$(3).
\end{proof}

\begin{lemma}[A capacity-dominated measure]\label{lem:tracemeasure}
Let $K\subset\mathbb R^n$ be compact with $\rev{\operatorname{Cap}_{p,1}}(K)>0$. There is a finite positive Radon measure $\nu$ supported on $K$ such that $\nu(K)>0$ and
\begin{equation}\label{eq:measurecap}
 \nu(A)\le\rev{\operatorname{Cap}_{p,1}}(A)\qquad\text{for every Borel set }A\subset\mathbb R^n.
\end{equation}
\end{lemma}
\begin{proof}
Let $h$ minimize the full Sobolev energy subject to $h\ge1$ quasieverywhere on $K$. Then $0\le h\le1$, and the distribution $-\Delta_ph+h^{p-1}$ is a positive measure $\nu$ supported on $K$. Testing with $h$ and using $h=1$ $\nu$-almost everywhere gives $\nu(K)=\rev{\operatorname{Cap}_{p,1}}(K)>0$. The energy-dual estimate shows that $\nu$ charges no set of zero $p$-capacity. For a nonnegative smooth $\psi$ of compact support, test the equilibrium equation with $\psi^p/(h+\varepsilon)^{p-1}$. Young's inequality gives
\[
 \int_{\mathbb R^n}\frac{\psi^p}{(h+\varepsilon)^{p-1}}\,d\nu
 \le\int_{\mathbb R^n}(|\nabla\psi|^p+|\psi|^p)\,dx.
\]
\rev{Sobolev approximation and truncation justify the test. Let $\varepsilon\downarrow0$ and use $h\le1$. Taking the infimum over capacitary cutoffs gives \eqref{eq:measurecap} on compact sets, hence on Borel sets by inner regularity of $\nu$.}
\end{proof}

\begin{proof}[Proof of (3)$\Rightarrow$(2) in Theorem~\ref{thm:trace}]
Suppose (2) fails. Quasicontinuity and capacitability give $\beta>0$ and a compact set $K\subset F$ of positive Sobolev capacity on which
$|\widetilde{\mathcal E_Du}|\ge2\beta$ outside a set of zero capacity. Let $\nu$ be the measure in Lemma~\ref{lem:tracemeasure}. By \eqref{eq:densitytrace}, \eqref{eq:lebesguetrace}, and H\"older's inequality, for $\nu$-almost every $y\in K$ and all sufficiently small $r$,
\begin{equation}\label{eq:tracepositive}
 \int_{D\cap B(y,r)}|u|^p\,dx\ge c \beta^pr^n.
\end{equation}
\rev{Indeed, the set on which $|u|<\beta$ inside $D\cap B(y,r)$ has measure $o(r^n)$ by \eqref{eq:lebesguetrace}. The density estimate \eqref{eq:densitytrace} then gives \eqref{eq:tracepositive}.}
Choose an integer $k_0\ge0$ such that \eqref{eq:tracepositive} holds for every $\rho_k=2^{-k}$, $k\ge k_0$, on a set of positive $\nu$-measure. By inner regularity, choose a compact subset $K'\subset K$ of this set with $\nu(K')>0$.

For $|x-y|<r$, the inclusion $B(y,r)\subset B(x,2r)$ and the comparison $\operatorname{Cap}_{p,r}\le2^p\operatorname{Cap}_{p,2r}$ give $\Theta_{F,p}(y,r)\le C_{n,p}\Theta_{F,p}(x,2r)$ by \eqref{eq:ballcap}. For $r\le1$, we also have $\operatorname{Cap}_{p,1}\le\operatorname{Cap}_{p,r}$, so Lemma~\ref{lem:tracemeasure} and the same capacity comparison give
\[
 \nu(K'\cap B(x,r))\le C_{n,p}r^{n-p}\Theta_{F,p}(x,2r).
\]

For $r=\rho_k$ with $k\ge k_0$, multiply \eqref{eq:tracepositive} by $r^{-n}\Theta_{F,p}(y,r)^{1/(p-1)}$, integrate over $K'$, and apply Tonelli's theorem. These bounds yield
\begin{equation}\label{eq:traceonescale}
 \int_{K'}\Theta_{F,p}(y,r)^{1/(p-1)}\,d\nu(y)
 \le\frac C{\beta^p}\int_D|u(x)|^p r^{-p}\Theta_{F,p}(x,2r)^{p/(p-1)}\,dx.
\end{equation}

By the fine Kellogg property \cite[Corollary~1.3]{BBL}, the set $F$ is not $p$-thin at $p$-quasievery point of $F$. The Wiener definition of thinness \cite[Definition~6.1]{BBL}, the capacity comparison \eqref{eq:relative-capacity}, and dyadic discretization imply
\begin{equation}\label{eq:wienertrace}
 \sum_{k\ge k_0}\Theta_{F,p}(y,\rho_k)^{1/(p-1)}=\infty
 \quad\text{for }\nu\text{-almost every }y\in K'.
\end{equation}
\rev{For $2r\le \rho<4r$, we have $B(x,2r)\subset B(x,\rho)$ and $\operatorname{Cap}_{p,\rho}\ge(2r/\rho)^p\operatorname{Cap}_{p,2r}$. Hence \eqref{eq:ballcap} gives $\Theta_{F,p}(x,\rho)\ge c_{n,p}\Theta_{F,p}(x,2r)$.} The intervals $[2\rho_k,4\rho_k)$ have disjoint interiors, and therefore
\begin{equation}\label{eq:tracesum}
 \sum_{k\ge k_0}\rho_k^{-p}\Theta_{F,p}(x,2\rho_k)^{p/(p-1)}\le C_{n,p}\rev{\bm{V}}_{F,p}(x).
\end{equation}
Integrate \eqref{eq:wienertrace} over $K'$. Tonelli's theorem, \eqref{eq:traceonescale}, and \eqref{eq:tracesum} give
\[
 \infty=\sum_{k\ge k_0}\int_{K'}\Theta_{F,p}(y,\rho_k)^{1/(p-1)}\,d\nu(y)
 \le\frac C{\beta^p}\int_D|u|^p\rev{\bm{V}}_{F,p}\,dx.
\]
This contradicts (3), completing the proof of Theorem~\ref{thm:trace}.
\end{proof}

\begin{corollary}\label{cor:fulltrace}
On a bounded $W^{1,p}$-extension domain,
\[
 \begin{aligned}
 W^{1,p}_0(D)&=\{u\in W^{1,p}(D):|u|^p\rev{\bm{V}}_{\mathbb R^n\setminus D,p}\in L^1(D)\}\\
 &=\{u\in W^{1,p}(D):|u|^p\rev{\bm{W}}_{D,p}\in L^1(D)\}.
 \end{aligned}
\]
\end{corollary}
\begin{proof}
The weighted bounds follow from Theorem~\ref{thm:hardy}. Conversely, use Theorem~\ref{thm:trace} with $F=\partial D$, the identity $W_{\partial D}^{1,p}(D)=W^{1,p}_0(D)$, and $\rev{\bm{V}}_{\partial D,p}\le \rev{\bm{V}}_{\mathbb R^n\setminus D,p}\le C_{n,p}\rev{\bm{W}}_{D,p}$.
\end{proof}
\endgroup

\begingroup\color{black}
\begin{corollary}[Distance to the Dirichlet portion]\label{cor:dirichlet-distance}
Under the assumptions of Theorem~\ref{thm:trace}, let $F\ne\varnothing$ and suppose that $F$ satisfies the $p$-capacity density condition up to a scale $0<r_0\le\operatorname{diam}D$. Thus, for $y\in F$ and $0<s<r_0$,
\[
 \operatorname{cap}_p(F\cap B(y,s),B(y,2s))
 \ge \eta\operatorname{cap}_p(B(y,s),B(y,2s))
\]
for some $\eta>0$. Write $d_F(x)=\operatorname{dist}(x,F)$. Then
\[
 \int_D\frac{|u(x)|^p}{d_F(x)^p}\,dx
 \le C\left(\int_D|\nabla u|^p\,dx+r_0^{-p}\int_D|u|^p\,dx\right),
 \qquad u\in W_F^{1,p}(D),
\]
where $C$ depends only on $n$, $p$, $\eta$, $D$, and the norm of the extension operator. Moreover, for every $u\in W^{1,p}(D)$,
\[
 u\in W_F^{1,p}(D)\quad\Longleftrightarrow\quad
 \int_D\frac{|u(x)|^p}{d_F(x)^p}\,dx<\infty.
\]
\end{corollary}
\begin{proof}
If $d_F(x)<r_0/2$, choose a nearest point $y\in F$. For $2d_F(x)\le r\le4d_F(x)$, we have $B(y,r/2)\subset B(x,r)$ and $r/2<r_0$. By \eqref{eq:relative-capacity}, \eqref{eq:ballcap}, and $\operatorname{Cap}_{p,r}\ge 2^{-p}\operatorname{Cap}_{p,r/2}$, the density bound at $y$ gives $\Theta_{F,p}(x,r)\ge c_{n,p}\eta$. Integrating over these radii gives $\bm{V}_{F,p}(x)\ge c_{n,p,\eta}d_F(x)^{-p}$. On the rest of $D$, $d_F^{-p}\le (2/r_0)^p$. The estimate follows from \eqref{eq:partialbound} and $r_0\le\operatorname{diam}D$. For the converse implication in the characterization, $\bm{V}_{F,p}\le d_F^{-p}$ by \eqref{eq:theta}. Thus finiteness of the distance-weight integral for $u\in W^{1,p}(D)$ implies $u\in W_F^{1,p}(D)$ by Theorem~\ref{thm:trace}.
The distance-weight estimate for $p=2$ recovers \cite[(8.13)]{Hara}.
\end{proof}

\endgroup

\section[Comparison of the Hardy weights]{\texorpdfstring{\rev{Comparison of the Hardy weights}}{Comparison of the Hardy weights}}\label{sec:comparison}
\begingroup\color{black}
\begin{theorem}[Comparison of the weights]\label{thm:comparison}
Let $D\subsetneq\mathbb R^n$ be nonempty and open, and let $1<p<\infty$. Then
\begin{equation}\label{eq:weightcomparison}
 \bm{W}_{D,p}(x)\asymp_{n,p}\bm{V}_{\mathbb R^n\setminus D,p}(x)
\end{equation}
for every $x\in D$.
\end{theorem}

For $1<p\le n$, the proof of $\bm{W}_{D,p}\le C_{n,p}\bm{V}_{\mathbb R^n\setminus D,p}$ uses an exponentially weighted estimate for the equilibrium measure in $\mathbb R^n$.
Exponential localization is a classical device in potential theory with a zeroth-order term, as reflected in the exponentially decaying Bessel kernels underlying nonlinear capacity theory; see Hedberg and Wolff \cite{HedbergWolff} and Adams and Hedberg \cite[Section~1.2.4]{AdamsHedberg}. Here we use $e^{-a|x-y|/r}$ directly as a test weight for the nonlinear equilibrium equation with absorption, obtaining the measure estimate \eqref{eq:expmeasure}.

\begin{lemma}\label{lem:exponentialmeasure}
Let $1<p<\infty$, let $K\subset\mathbb R^n$ be compact, and let $h_r$ minimize the energy defining $\operatorname{Cap}_{p,r}(K)$ over functions in $W^{1,p}(\mathbb R^n)$ that are at least $1$ $p$-quasieverywhere on $K$. Write
\[
 \mu_r=-\Delta_p h_r+r^{-p}h_r^{p-1}.
\]
There is $a=a(n,p)\in(0,1]$ such that, for every $x\in\mathbb R^n$ and $r>0$,
\begin{equation}\label{eq:expmeasure}
 r^{-p}\int_{\mathbb R^n}e^{-a|x-y|/r}h_r(y)^{p-1}\,dy
 \le2\int_{\mathbb R^n}e^{-a|x-y|/r}\,d\mu_r(y).
\end{equation}
Moreover, $0\le h_r\le1$, the measure $\mu_r$ is nonnegative and supported on $K$, and
\begin{equation}\label{eq:equilibriumcap}
 \mu_r(A)\le\operatorname{Cap}_{p,r}(K\cap A)
\end{equation}
for every Borel set $A\subset\mathbb R^n$.
\end{lemma}
\begin{proof}
Truncation gives $0\le h_r\le1$. The variational inequality shows that $\mu_r$ is a nonnegative Radon measure supported on $K$. The argument of Lemma~\ref{lem:tracemeasure}, applied at scale $r$, gives \eqref{eq:equilibriumcap}.

We first show that, for every ball $B=B(z,r)$,
\begin{equation}\label{eq:gradmass}
 \int_B|\nabla h_r|^{p-1}\,dy
 \le C_{n,p}r^{1-p}\int_{2B}h_r^{p-1}\,dy.
\end{equation}
Fix a sufficiently small $\eta\in(0,p-1)$ such that the weak Harnack inequality for nonnegative supersolutions holds with exponent $\gamma=(1+\eta)(p-1)>p-1$. Put $v=h_r+\varepsilon$, where $\varepsilon>0$. Since $-\Delta_p v+r^{-p}v^{p-1}\ge0$, we may test the supersolution inequality with $\chi^pv^{-\eta}$, where $\chi\in C_c^\infty(2B)$, $0\le\chi\le1$, $\chi=1$ on $B$, and $|\nabla\chi|\le C/r$. Young's inequality gives
\[
 \int_B|\nabla h_r|^pv^{-1-\eta}\,dy
 \le Cr^{-p}\int_{2B}v^{p-1-\eta}\,dy.
\]
By H\"older's inequality,
\[
 \int_B|\nabla h_r|^{p-1}\,dy
 \le Cr^{1-p}\left(\int_{2B}v^{p-1-\eta}\,dy\right)^{(p-1)/p}
 \left(\int_Bv^{\gamma}\,dy\right)^{1/p}.
\]
The weak Harnack inequality \cite[Theorem~1.2]{Trudinger} and H\"older's inequality, applied on concentric balls at scale $r$, yield
\begin{align*}
 \left(\frac1{|B|}\int_Bv^\gamma\,dy\right)^{1/\gamma}
 &\le C\left(\frac1{|2B|}\int_{2B}v^{p-1}\,dy\right)^{1/(p-1)},\\
 \frac1{|2B|}\int_{2B}v^{p-1-\eta}\,dy
 &\le\left(\frac1{|2B|}\int_{2B}v^{p-1}\,dy\right)^{(p-1-\eta)/(p-1)}.
\end{align*}
Substitution, followed by $\varepsilon\downarrow0$, proves \eqref{eq:gradmass}. The constants are uniform because the coefficient of the absorption term is $r^{-p}$.

Fix $x$ and put $\omega(y)=e^{-a|x-y|/r}$. If $a\le1$, the values of $\omega$ on $2B$ are comparable to $\omega(z)$, with a uniform constant. Multiply \eqref{eq:gradmass} by $\omega(z)$, integrate with respect to $z$, and use Tonelli's theorem to obtain
\begin{equation}\label{eq:weightedgradient}
 \int_{\mathbb R^n}\omega|\nabla h_r|^{p-1}\,dy
 \le Cr^{1-p}\int_{\mathbb R^n}\omega h_r^{p-1}\,dy.
\end{equation}
Test the equation for $h_r$ with a compactly supported approximation of $\omega$. The cutoff error tends to zero, since $h_r\in W^{1,p}(\mathbb R^n)$, $\mu_r$ is finite, and \eqref{eq:weightedgradient} holds. Using $|\nabla\omega|\le(a/r)\omega$ gives
\begin{align*}
 r^{-p}\int_{\mathbb R^n}\omega h_r^{p-1}\,dy
 &\le\int_{\mathbb R^n}\omega\,d\mu_r
 +\frac ar\int_{\mathbb R^n}\omega|\nabla h_r|^{p-1}\,dy\\
 &\le\int_{\mathbb R^n}\omega\,d\mu_r
 +Ca r^{-p}\int_{\mathbb R^n}\omega h_r^{p-1}\,dy.
\end{align*}
Choose $a=a(n,p)>0$ so that $Ca\le1/2$. Absorbing the last term proves \eqref{eq:expmeasure}.
\end{proof}

\begin{lemma}\label{lem:equilibriumscales}
Suppose $1<p\le n$, and let $h_r$ and $\mu_r$ be as in Lemma~\ref{lem:exponentialmeasure}. For every $x\notin K$ and $r>0$,
\begin{equation}\label{eq:equilibriumscales}
 \begin{aligned}
 h_r(x)\le C_{n,p}\bigg[&\int_0^r
 \left(\frac{\mu_r(B(x,s))}{s^{n-p}}\right)^{1/(p-1)}\frac{ds}{s}\\
 &+\left(r^{p-n}\int_{\mathbb R^n}e^{-a|x-y|/r}\,d\mu_r(y)\right)^{1/(p-1)}\bigg],
 \end{aligned}
\end{equation}
where $a$ is the constant in \eqref{eq:expmeasure}.
\end{lemma}
\begin{proof}
Apply the interior potential estimate \cite[Theorem~3.3]{Maly96} to
$-\Delta_ph_r+r^{-p}h_r^{p-1}=\mu_r$ in $\mathbb R^n$, with outer radius $r$. Rescaling to $r=1$ makes the structural constants independent of $r$. Since $x\notin K$, the solution is continuous near $x$. Thus, for a fixed $\gamma\in(p-1,n(p-1)/(n-p+1))$,
\[
 h_r(x)\le C\left[\left(\frac1{|B(x,r)|}\int_{B(x,r)}h_r^\gamma\,dy\right)^{1/\gamma}
 +\int_0^r\left(\frac{\mu_r(B(x,s))}{s^{n-p}}\right)^{1/(p-1)}\frac{ds}{s}\right].
\]
The boundary-capacity term in that theorem vanishes because the equation is posed on all of $\mathbb R^n$. Choose $\gamma$ also within the weak Harnack range. As $h_r$ is a nonnegative supersolution of the homogeneous absorption equation,
\[
 \left(\frac1{|B(x,r)|}\int_{B(x,r)}h_r^\gamma\,dy\right)^{1/\gamma}
 \le C\left(\frac1{|B(x,2r)|}\int_{B(x,2r)}h_r^{p-1}\,dy\right)^{1/(p-1)}.
\]
On $B(x,2r)$, the weight in \eqref{eq:expmeasure} is at least $e^{-2a}$. Therefore \eqref{eq:expmeasure} bounds the last expression by
\[
 C\left(r^{p-n}\int_{\mathbb R^n}e^{-a|x-y|/r}\,d\mu_r(y)\right)^{1/(p-1)}.
\]
These estimates prove \eqref{eq:equilibriumscales}.
\end{proof}

\begin{lemma}[Comparison with a ball]\label{lem:ballresponse}
For every nonempty proper open set $D\subset\mathbb R^n$, $x\in D$, and $r>0$, we have
\begin{equation}\label{eq:ballbound}
 z_{r,p}^D(x)\le
 \left[1+\frac{p-1}{p}
 \left(\frac{\delta_D(x)^p}{nr^p}\right)^{1/(p-1)}\right]^{-1}.
\end{equation}
\end{lemma}
\begin{proof}
By domain comparison, $z_{r,p}^D(x)\le z_{r,p}^{B(x,\delta_D(x))}(x)$. Write $z(\rho)$ for the radial solution in $B(x,\delta_D(x))$, measured from its center. The radial equation and monotonicity give
\[
 (z'(\rho))^{p-1}=r^{-p}\rho^{1-n}
 \int_0^\rho t^{n-1}z(t)^{p-1}\,dt
 \ge\frac{r^{-p}\rho}{n}z(0)^{p-1}.
\]
Integrating from $0$ to $\delta_D(x)$ and using $z(\delta_D(x))=1$ proves \eqref{eq:ballbound}.
\end{proof}

\begin{proof}[Proof of Theorem~\ref{thm:comparison}]
Theorem~\ref{thm:hardy} gives $\bm{V}_{\mathbb R^n\setminus D,p}\le C_{n,p}\bm{W}_{D,p}$. We prove the reverse bound first for $1<p\le n$ and bounded $D$.

Choose $R$ with $\overline D\subset B(0,R)$ and set $K=\overline{B(0,R)}\setminus D$. Let $h_r$ be as in Lemma~\ref{lem:exponentialmeasure}. Multiplying $1-h_r$ by a smooth cutoff equal to $1$ near $\overline D$ gives a function in $W^{1,p}(\mathbb R^n)$ that vanishes $p$-quasieverywhere outside $D$. The zero-boundary characterization \cite[Theorem~4.5]{HKM} and uniqueness in \eqref{eq:absorption} give $h_r|_D=z_{r,p}^D$.

Fix $x\in D$ and set
\[
 f(s)=s^{-1}\Theta_{K,p}(x,s)^{1/(p-1)}.
\]
For $0<s\le r$, \eqref{eq:equilibriumcap} and the defining energies show that
\begin{align*}
 \mu_r(B(x,s))&\le\operatorname{Cap}_{p,r}(K\cap B(x,s))\\
 &\le\operatorname{Cap}_{p,s}(K\cap B(x,s))
 =C_{n,p}s^{n-p}\Theta_{K,p}(x,s).
\end{align*}
To estimate the second term of \eqref{eq:equilibriumscales}, put $s_j=2^{j+1}r$ for $j\ge0$ and decompose $K$ into $K\cap B(x,2r)$ and the annuli $K\cap(B(x,s_j)\setminus B(x,2^jr))$, $j\ge1$. The measure domination \eqref{eq:equilibriumcap} and the comparison of energies give
\begin{align*}
 \mu_r(B(x,s_j))&\le\operatorname{Cap}_{p,r}(K\cap B(x,s_j))\\
 &\le(s_j/r)^p\operatorname{Cap}_{p,s_j}(K\cap B(x,s_j))
 =C_{n,p}r^{-p}s_j^n\Theta_{K,p}(x,s_j).
\end{align*}
Consequently,
\[
 r^{p-n}\int_{\mathbb R^n}e^{-a|x-y|/r}\,d\mu_r(y)
 \le C_{n,p}\sum_{j\ge0}e^{-a2^j}2^{jn}\Theta_{K,p}(x,s_j).
\]
Fix $q>0$. If $q\ge1$, apply H\"older's inequality; if $0<q\le1$, use the subadditivity of the $q$th power. In either case, after decreasing the exponential constant, we obtain
\[
 \left(\sum_{j\ge0}e^{-a2^j}2^{jn}b_j\right)^q
 \le C\sum_{j\ge0}e^{-c2^j}2^{jnq}b_j^q
 \qquad(b_j\ge0).
\]
Apply this with $q=1/(p-1)$ and absorb the remaining polynomial factor in the exponential. Dividing \eqref{eq:equilibriumscales} by $r$ now yields
\begin{equation}\label{eq:logconvolution}
 \frac{z_{r,p}^D(x)}r
 \le C_{n,p}\left[\int_0^r\frac sr f(s)\,\frac{ds}{s}
 +\sum_{j\ge0}2^{j+1}e^{-c2^j}f(2^{j+1}r)\right].
\end{equation}
The first operator on the right is convolution with $e^{-t}\bm{1}_{(0,\infty)}(t)$ in logarithmic coordinates. Its kernel has $L^1$ norm $1$. The coefficients of the second operator are summable, and every dilation $r\mapsto2^{j+1}r$ preserves $dr/r$. Young's convolution inequality and Minkowski's inequality therefore give
\begin{align*}
 \bm{W}_{D,p}(x)
 &=p\int_0^\infty\left(\frac{z_{r,p}^D(x)}r\right)^p\frac{dr}{r}\\
 &\le C_{n,p}p\int_0^\infty f(s)^p\,\frac{ds}{s}
 =C_{n,p}\bm{V}_{K,p}(x)
 \le C_{n,p}\bm{V}_{\mathbb R^n\setminus D,p}(x).
\end{align*}

For $p>n$, Lemma~\ref{lem:ballresponse} and the substitution $t=\delta_D(x)/r$ give
\[
 \begin{aligned}
 \bm{W}_{D,p}(x)
 &\le p\delta_D(x)^{-p}\int_0^\infty
 \frac{t^{p-1}}{\bigl(1+c_{n,p}t^{p/(p-1)}\bigr)^p}\,dt\\
 &\le C_{n,p}\delta_D(x)^{-p}.
 \end{aligned}
\]
Choose a nearest point $y\in \mathbb R^n\setminus D$ to $x$. Since $p>n$, a singleton has positive $\operatorname{Cap}_{p,1}$-capacity. For every $s>\delta_D(x)$, scaling and \eqref{eq:ballcap} give $\Theta_{\mathbb R^n\setminus D,p}(x,s)\ge c_{n,p}$. Integrating in \eqref{eq:theta} and using \eqref{eq:capupper} and Theorem~\ref{thm:hardy}, we obtain
\begin{equation}\label{eq:supercritical}
 \bm{V}_{\mathbb R^n\setminus D,p}(x)\asymp_{n,p}\bm{W}_{D,p}(x)\asymp_{n,p}\delta_D(x)^{-p}.
\end{equation}
In particular, the upper bound in \eqref{eq:weightcomparison} holds for $p>n$.

Finally, for unbounded $D$ and $1<p\le n$, set $D_j=D\cap B(0,j)$. For fixed $x\in D$ and $s>0$, we have $(\mathbb R^n\setminus D_j)\cap B(x,s)=(\mathbb R^n\setminus D)\cap B(x,s)$ for all sufficiently large $j$. Since $0\le\Theta_{\mathbb R^n\setminus D_j,p}\le1$ and these densities vanish below a fixed positive radius for all sufficiently large $j$, dominated convergence gives $\bm{V}_{\mathbb R^n\setminus D_j,p}(x)\to\bm{V}_{\mathbb R^n\setminus D,p}(x)$. By the definition of $z_{r,p}^D$ and Fatou's lemma,
\[
 \bm{W}_{D,p}(x)\le\liminf_{j\to\infty}\bm{W}_{D_j,p}(x)
 \le C_{n,p}\bm{V}_{\mathbb R^n\setminus D,p}(x).
\]
The bound for $p>n$ did not require boundedness of $D$. This completes the proof.
\end{proof}
\endgroup
\section{The limit as \texorpdfstring{$p$}{p} decreases to one}\label{sec:one}
\begingroup\color{black}
Let $D$ be a nonempty bounded Lipschitz domain. Write $P(A)$ for the perimeter in $\mathbb R^n$, $v^0$ for the zero extension of $v$ from $D$, and $\operatorname{TV}(v^0)$ for its total variation in $\mathbb R^n$. For $\lambda>0$, let $A_\lambda$ be the maximal minimizer, modulo null sets, of
\begin{equation}\label{eq:setmin}
 m(\lambda)=\min_{A\subset D}\{P(A)+\lambda|D\setminus A|\}.
\end{equation}
Applying Barozzi's construction \cite[Theorem~2.1 and Remark~2.4]{Barozzi} with density $1$ on $D$, we obtain a nested family and the following almost-everywhere formula for the interior variational mean curvature:
\begin{equation}\label{eq:oneweight}
 -\rev{\mathsf{H}}_D(x)=\inf\{\lambda>0:x\in A_\lambda\}
 =\int_0^\infty(1-\bm{1}_{A_\lambda}(x))\,d\lambda.
\end{equation}
One may use the nested family at positive rational parameters to define the threshold and choose measurable representatives. The function $m$ is concave. At each differentiability point, every minimizer $A$ has $|D\setminus A|=m'(\lambda)$, by comparison with its affine energy as a function of $\lambda$. Since $A_\lambda$ contains every minimizer, the minimizer is unique modulo null sets for almost every $\lambda$.

\begin{theorem}\label{thm:one}
As $p\downarrow1$,
\begin{equation}\label{eq:oneconv}
 \rev{\bm{W}}_{D,p}\longrightarrow -\rev{\mathsf{H}}_D\quad\text{in }L^1_{\mathrm{loc}}(D).
\end{equation}
\end{theorem}
\endgroup

\begin{lemma}\label{lem:fixedone}
Let $v_{p,\lambda}$ minimize
\[
 \int_D|\nabla v|^p\,dx+\lambda\int_D|1-v|^p\,dx,
 \qquad v\in W^{1,p}_0(D).
\]
For almost every $\lambda>0$,
$(1-v_{p,\lambda})^p\to1-\rev{\bm{1}}_{A_\lambda}$ in $L^1(D)$ as $p\downarrow1$.
\end{lemma}
\begin{proof}
Truncation gives $0\le v_{p,\lambda}\le1$, and comparison with zero bounds its energy by $\lambda|D|$. H\"older's inequality bounds the total variations of the zero extensions uniformly for $p$ near $1$. \rev{By BV compactness, the zero extensions have subsequences converging in $L^1(\mathbb R^n)$.} The inequality $t^p\ge pt-(p-1)$ and lower semicontinuity yield the liminf bound for the functional
\begin{equation}\label{eq:bvfunctional}
 \rev{\operatorname{TV}(v^0)}+\lambda\int_D(1-v)\,dx,\qquad 0\le v\le1.
\end{equation}
\rev{We use the boundary approximation property of a Lipschitz domain \cite[Chapter~3]{AFP}. Every such $v$ with finite zero-extension variation admits compactly supported Sobolev approximations in $D$ that converge in $L^1(D)$ and whose total variations converge to $\operatorname{TV}(v^0)$.} Smooth approximation followed by a diagonal choice as $p\downarrow1$ gives a recovery sequence for \eqref{eq:bvfunctional}. Thus every limit of minimizers minimizes this functional.

The coarea formula writes \eqref{eq:bvfunctional} as
\[
 \int_0^1\bigl[P(\{v>t\})+\lambda|D\setminus\{v>t\}|\bigr]\,dt.
\]
\rev{For almost every $\lambda$, uniqueness of $A_\lambda$ forces the superlevel sets of every minimizer to agree with $A_\lambda$ for almost every $t\in(0,1)$. Hence the unique minimizer of \eqref{eq:bvfunctional} is $\bm{1}_{A_\lambda}$.}
All subsequential limits coincide. Since $t^p\to t$ uniformly on $[0,1]$, the convergence in \rev{$L^1(D)$} follows.
\end{proof}

\begin{proof}[Proof of \eqref{eq:oneconv}]
\rev{Since $v_{p,r^{-p}}=J_{r,p}1$, the substitution $\lambda=r^{-p}$ gives}
\begin{equation}\label{eq:lambdarepresentation}
 \rev{\bm{W}}_{D,p}(x)=\int_0^\infty(1-v_{p,\lambda}(x))^p\,d\lambda.
\end{equation}
\begingroup\color{black}
For the tail of \eqref{eq:lambdarepresentation}, use the ball estimate \eqref{eq:ballbound} from Lemma~\ref{lem:ballresponse}.

Fix a compact $K\subset D$. The integral over $(0,\varepsilon)$ is at most $\varepsilon$ pointwise. Put $q=p-1$ and $\sigma=r/\delta_D(x)$. By \eqref{eq:ballbound}, for $0<q\le1/2$ and $n\sigma^{1/2}\le1/2$,
\[
 z_{r,p}^D(x)\le \frac{p}{q}(n\sigma^{1+q})^{1/q}
 =\frac{p}{q}\sigma^2(n\sigma^{1-q})^{1/q}\le C \sigma^2.
\]
Since $\inf_K\delta_D>0$, this gives $(1-v_{p,\lambda})^p\le C_K\lambda^{-2}$ on $K$ for large $\lambda$, uniformly in $p$ near $1$. Thus Lemma~\ref{lem:fixedone} and dominated convergence in $L^1(K)$ give the limit in \eqref{eq:oneweight}, proving \eqref{eq:oneconv}.
\endgroup
\end{proof}

\begingroup\color{black}
\begin{corollary}[Sharp BV Hardy inequality]\label{cor:bvhardy}
For every $u\in BV(D)$,
\begin{equation}\label{eq:bvhardy}
 \int_D(-\rev{\mathsf{H}}_D)|u|\,dx\le \operatorname{TV}(u^0).
\end{equation}
The constant $1$ is optimal.
\end{corollary}
\begin{proof}
Barozzi's mass formula and variational inequality give
\begin{equation}\label{eq:massone}
 \int_D(-\rev{\mathsf{H}}_D)\,dx=P(D),\qquad
 \int_A(-\rev{\mathsf{H}}_D)\,dx\le P(A)\quad(A\subset D).
\end{equation}
Indeed, $D$ minimizes $P(A)-\int_A(-\rev{\mathsf{H}}_D)\,dx$ among $A\subset D$, and its value at $D$ is zero \cite[Theorem~2.1]{Barozzi}. Applying \eqref{eq:massone} to the superlevel sets of $|u^0|$ and using the coarea formula yields \eqref{eq:bvhardy}.
The constant $1$ is optimal by taking $u=1$ on $D$.
\end{proof}
The total variation in \eqref{eq:bvhardy} includes the boundary trace term \cite[Chapter~3]{AFP}. In particular, $u=1$ has finite weighted integral and nonzero trace, so the trace characterization in Theorem~\ref{thm:trace} fails at the endpoint.

On a ball of radius $R$, the isoperimetric inequality gives $-\rev{\mathsf{H}}_D=n/R$. On the same ball the complement contains a ball of radius $\delta_D(x)$ within $B(x,4\delta_D(x))$ for every $x\in D$. By scaling and monotonicity, $\Theta_{\mathbb R^n\setminus D,p}(x,r)\ge c_{n,p}>0$ for $4\delta_D(x)\le r\le8\delta_D(x)$. Thus $\rev{\bm{W}}_{D,p}\ge c_{n,p}\delta_D^{-p}$ by Theorem~\ref{thm:hardy}, and $\rev{\bm{W}}_{D,p}\notin L^1(D)$ for $p>1$. Consequently, the convergence in \eqref{eq:oneconv} cannot, in general, be strengthened to convergence in $L^1(D)$.
\endgroup

\backmatter
% Keep acknowledgements as a top-level PDF bookmark; \bmhead is a run-in
% subparagraph heading in the Springer Nature class.
\pdfbookmark[1]{Acknowledgements}{acknowledgements}
\hypersetup{bookmarksdepth=1}
\bmhead{Acknowledgements}
This work was supported by JSPS KAKENHI Grant Numbers 21K03281 and 23H03798.
The author used ChatGPT (OpenAI) to assist with mathematical discussions,
manuscript preparation, and English-language editing.

\pdfbookmark[1]{Data availability}{data-availability}

\bmhead{Conflict of Interest}
The author declares that there is no conflict of interest.

\bmhead{Data availability}
No datasets were generated or analysed during the current study.

\pdfbookmark[1]{References}{references}

\end{document}